\documentclass[a4paper,12pt]{article}
\usepackage[utf8]{inputenc}
\usepackage[T1]{fontenc}
\usepackage{amsmath,amsfonts,amssymb,amsthm,amscd,tipa,color,graphicx,tikz-cd,longtable}
\usepackage[english]{babel}
\usepackage[a4paper,margin=0.75in,top=4.5cm,bottom=3.5cm]{geometry}

\renewcommand{\P}{\mathbb{P}}
\renewcommand{\S}{\mathbb{S}}
\newcommand{\C}{\mathbb{C}}
\newcommand{\hc}{\textbf H^2_\C}

\DeclareMathOperator{\aut}{Aut}

\DeclareMathOperator{\tgt}{T}
\DeclareMathOperator{\is}{Isom}
\DeclareMathOperator{\Mod}{Mod}
 \DeclareMathOperator{\Isom}{Isom}

\theoremstyle{plain}
\begingroup

\newtheorem{theo}{Theorem}[section]
\newtheorem*{theo*}{Theorem}
\newtheorem{corollary}[theo]{Corollary}

\newtheorem{lemma}[theo]{Lemma}
\newtheorem{defin}[theo]{Definition}
\newtheorem{rk}[theo]{Remark}
\newtheorem{prop}[theo]{Proposition}

\title{Complex hyperbolic orbifolds and Lefschetz fibrations}
\author{Elisha Falbel and Irene Pasquinelli}
\date{\vspace{-3ex}}

\begin{document}

\maketitle

\begin{abstract}
A class of complex hyperbolic lattices in $PU(2,1)$ called the Deligne-Mostow lattices has been reinterpreted by Hirzebruch (see \cite{hirz}, \cite{BHH} and \cite{tretkoff}) in terms of line arrangements.
They use branched covers over a suitable blow up of the complete quadrilateral arrangement of lines in $\P^2$ to construct the complex hyperbolic surfaces over the orbifolds associated to the lattices. 

In \cite{irene3F} and \cite{irene2}, fundamental domains for these lattices have been built by Pasquinelli. 
Here we show how the fundamental domains can be interpreted in terms of line arrangements as above.  
This parallel is then applied in two contexts. 

Dashyan in \cite{ruben} uses Hirzebruch's construction to build infinitely many representations of 3-manifolds.
Here we show that his construction can be generalised to all of the Deligne-Mostow lattices and more representations can be built. 

Wells in \cite{wells} shows that two of the Deligne-Mostow lattices in $PU(2,1)$ can be seen as hybrids of lattices in $PU(1,1)$. 
Here we show that he implicitly uses the line arrangement and we complete his analysis to all possible pairs of lines. 
In this way, we show that three more Deligne-Mostow lattices can be given as hybrids.

\end{abstract}

\tableofcontents

\section{Introduction}

The \emph{complex hyperbolic space} $\textbf H^n_\C$ is the complex equivalent to the real hyperbolic space.
It is defined starting from a Hermitian form of signature $(n,1)$ on $\C^{n+1}$ and it contains the projectivised of points of $\C^{n+1}$ which are negative for the product induced by the Hermitian form. 
Its isometry group is generated by the complex conjugation and the group $PU(n,1)$ of the projective matrices which preserve the Hermitian form. 

An important problem in complex hyperbolic geometry is the study of \emph{lattices} in $PU(n,1)$. 
To produce lattices, it is quite standard to use number theory to construct \emph{arithmetic} lattices.
% (i.e. lattices that are discrete in $PU(n,1)$ ``in the same'' way as $\Z$ is discrete in $\R$). 
Finding non-arithmetic lattices is an important and difficult problem . 

In the real case, Gromov and Piatetski-Shapiro in \cite{nonarithm} constructed infinitely many non-arithmetic lattices in any dimension, using a construction called \emph{hybridisation}. 
It consists in getting a new non-arithmetic lattice by ``gluing'' two arithmetic lattices along a totally geodesic hypersurface (a maximal dimension subspace which contains whole geodesics). 
In the complex case, an equivalent construction is not possible and the existence of non-arithmetic complex hyperbolic lattices is still an open question for $n \geq4$ and even in lower dimension very few constructions are known. 

%The complex hyperbolic space is an example of symmetric space (riemannian manifolds which admit at each point in isometry which is an involution and for which the point is a fixed point). 
%Symmetric spaces are classifies according to the rank (the maximal dimension of a flat subspace) and the type (the curvature). 
%Here we will be interested only in the spaces of non-compact type (with negative curvature). 
%We know that every arithmetic subgroup is a lattice.  
%Viceversa, if the rank is 2 or more, then all lattices are arithmetic. 
%The symmetric spaces of non compact type and of rank 1 are the four hyperbolic spaces (real hyperbolic space, complex hyperbolic space, quaternionic hyperbolic space and octonionic hyperbolic plane). 
%Among these, the complex hyperbolic space is the only case where the arithmeticity of lattices is not completely known. 
%In fact, for the last two cases, lattices are always arithmetic.

 The first examples of lattices were given by Picard using monodromy of hypergeometric functions.  
Those lattices were studied, among others, by Deligne and Mostow (see, for example, \cite{mostow}, \cite{delignemostow}).  
We refer to the list in \cite{delignemostow} as \emph{Deligne-Mostow lattices}.
In \cite{thurston}, Thurston gave a geometric interpretation of the same lattices. 

He started with cone metrics on a sphere with $N$ prescribed cone angles and area one. 
He proved that the area is a Hermitian form of signature $(N-3,1)$ and that the metric completion of the moduli space of such cone metrics has a complex hyperbolic structure of finite volume. 
Using authomorphisms of the sphere he then got an explicit list of cone singularities for which one can get a lattice. 
In this way, he produces the same lattices as the Deligne-Mostow construction. 
In \cite{irene3F} and \cite{irene2}, Pasquinelli generalised a construction introduced by Parker in \cite{livne} to build explicit fundamental domain for all of the Deligne-Mostow lattices in $PU(2,1)$ with 3- and 2-fold symmetry respectively.
This means that 3 (resp. 2) of the five cone points have the same cone angle.
All of the Deligne-Mostow lattices in dimension 2 have either 2- or 3-fold symmetry.

These lattices also have a third interpretation, introduced by Hirzebruch (see, for example, \cite{hirz}) for one case and generalised in \cite{BHH} (see also \cite{tretkoff}) for all of them. 
This construction is explained in Section \ref{sec:covers}.

Hirzebruch considers four points in $\P^2$ such that no three of them lie on the same line. 
This configuration determines six lines passing through pairs of those four points. 
In this way one defines an arrangement of six lines in $\P^2$ with four points of triple intersection and three points of double intersection (see Figure \ref{fig:arrang}). 
Blowing up the four points of triple intersection in $\P^2$, one obtains $\widehat{\P^2}$, with ten exceptional lines. 
One then defines a branched cover of $\P^2$, ramified along the ten lines with a fixed ramification index. 
A computation of the first and second Chern classes $c_1$ and $c_2$ of the branched cover shows that 
certain branched covers are complex hyperbolic surfaces.  
This is done using the equality case in Bogomolov-Miyaoka-Yau inequality, which guarantees that a compact complex surface \emph{of general type} (which is our case), satisfying $c_1^2-3c_2=0$, is complex hyperbolic. 
In fact, one proves that (\cite{BHH}) the surface is the quotient of (a torsion-free subgroup of) a Deligne-Mostow lattice.

The space $\widehat{\P^2}$ corresponds to the orbifold obtained by gluing the sides of the polyhedra which are the fundamental domains built in \cite{irene2}. 
This construction and the fact that the fundamental domains are very explicit allows us to see them in the blow up of $\P^2$. 
In Section \ref{fund} we built a precise ``dictionary'' between the line arrangement in $\widehat{\P^2}$ and the polyhedra. 
In particular one has

\begin{theo}Given any Deligne-Mostow lattice, there exists a fundamental domain such that (up to quotienting by symmetries) its 0-skeleton corresponds to the intersection points of the arrangement of lines in $\widehat{\P^2}$, and such that each segment of the 1-skeleton and each 2-facet which is a triangle projects to a line in the arrangement.
\end{theo}

This parallel can be applied to two different contexts to deduce some interesting results. 

The first one is the one of spherical CR-structures, which is a very rich topic related to complex hyperbolic geometry. 
We say that a 3-manifold $M$ has a spherical CR-structure if it is modelled on $\S^3$ seen as the boundary of $\hc$. 
Then one can consider the holonomy representation $\rho \colon \pi_1(M) \to \is(\hc)$. 
An important problem is whether the CR-structure is uniformisable, i.e. if one can write $M=\rho(\pi_1(M)) \textbackslash \Omega$, for some open $\Omega \subset \S^3$.
Very few examples of such representations (uniformisable or not) are known. 

In \cite{ruben}, Dashyan considers Hirzebruch's quadrilateral arrangement of lines.
He then shows that one can construct a Lefschetz fibration using the pencil of conics through the four points. 
He studies the monodromy of the fibration over a non-singular curve and considers the surface bundle over the circle with fibre the general fibre above. 
Using this, he builds infinitely many representations of the 3-manifold fundamental group into $PU(2,1)$ and studies the limit set.

In Section \ref{sec:reps}, we describe the fibration inside the fundamental domains from \cite{irene2} and we show that the work from \cite{ruben} can be generalised to all of the Deligne-Mostow lattices and one can, in this way, construct other representations. 

In particular, we obtain the following

\begin{theo}
For each element of $\Mod_{0,4}$, consider a surface bundle $M$ with the monodromy as above and with fibre the orbifold with the sphere as underlying space and with isotropy of order equal to the multiplicity of the exceptional fibre at each of the four marked points. Then there exists a representation of the orbifold fundamental group of $M$ into a Deligne-Mostow lattice in $\Isom(\hc)$.
\end{theo}

Here $\Mod_{0,4}$ is the mapping class group of a sphere with four marked points and $M$ is 
a 3-orbifold fibred over the circle obtained by gluing the sphere orbifold through an element of $\Mod_{0,4}$.

The second application is explained in Section \ref{sec:hybrid}. 
It is about the study of a possible complex equivalent to the hybridisation of Gromov and Piatetski-Shapiro \cite{nonarithm}. 
In the complex hyperbolic case, such a construction is impossible, since there are no totally geodesic real hypersurfaces along which to glue the lattices. 

In \cite{wells}, Wells modifies a possible generalisation to the complex case given by Hunt and already used by Paupert and shows that two of the Deligne-Mostow lattices in $PU(2,1)$ with 3-fold symmetry and of a certain type can be written as hybrids of two lattices in $PU(1,1)$, arithmetic and non-commensurable. 
One of the difficulties in his work is to choose a pair of totally geodesic hypersurfaces. 
 
Using line arrangements, we obtain natural candidates along which to hybridise. 
In this way, we show that three more (commensurability classes of) Deligne-Mostow lattices in $PU(2,1)$ can be described as hybrids.  
In fact, other than the two cases $(4,6)$ and $(5,4)$ found by Wells, also the lattices $(3,8)$, $(6,4)$ and $(3,4,4)$ are hybrids of two non-commensurable arithmetic lattices in $PU(1,1)$. 

We considered all lines in the arrangement of 10 exceptional curves in $\widehat{\P^2}$
and consider all 15 pairwise orthogonal intersections of these lines.  
The stabilizer of each line $L$ is a triangle group $\Delta$ and we obtained the following 

\begin{theo}
Let $\Gamma\subset PU(2,1)$ be a non-arithmetic Deligne-Mostow lattice.
Let $(L_1, L_2)$ be an orthogonally intersecting pair of lines with corresponding triangle groups $\Delta_1$ and $\Delta_2$.  
Then 
\begin{enumerate}
\item 
 either both $\Delta_1$ and $\Delta_2$ are arithmetic and commensurable,
\item or both $\Delta_1$ and $\Delta_2$ are arithmetic and non-commensurable and in this case $\Gamma$ is an hybridisation of $\Delta_1$ and $\Delta_2$,
 \item or one of the triangle groups is non-arithmetic.
\end{enumerate}

\end{theo}

\section{Deligne-Mostow lattices as branched covers of line arrangements}\label{sec:covers}

In this section we will explain Hirzebruch's construction (for a complete introduction to the subject see
 \cite{BHH} and \cite{tretkoff}. 
We consider an arrangement of lines in $\P^2$ and construct branched covers. 
For certain ramification orders these covers are complex hyperbolic ball quotients and we will explain how they are related to the Deligne-Mostow lattices.

\subsection{Branched covers of line arrangements}\label{sec:tretkoff}

For Hirzebruch's construction we will follow Dashyan \cite{ruben}, where one can also find detailed proofs. 
The goal is to build smooth complex algebraic surfaces using branched covering spaces. 

For an integer $n$, consider the morphism 
\begin{align*}
c_n \colon \P^{k-1} &\to \P^{k-1} \\
[u_1 \colon \dots \colon u_k] &\mapsto [u_1^n\colon \dots \colon u_k^{n}].
\end{align*}
This is a branched covering map of degree $n^{k-1}$.

The ramification locus is the arrangement of $k$ hyperplanes $D_i$ defined by $u_1=0$.
The preimage of the ramification locus is called the branching locus.
A point outside of the ramification locus has $n^{k-1}$ preimages. 
The number of preimages of a point contained in the ramification locus is $n^{k-1-r}$, if the point belongs to a certain number $r$ of hyperplanes $D_i$. 

We now want to create a branched covering of $\P^2$, ramified along a prescribed arrangement of $k$ lines $L_1=0, \dots, L_k=0$, with a similar behaviour to the one of the branching map $c_n$. 
We will consider the map $L \colon \P^2 \to \P^{k-1}$ sending $[z] \in \P^2$ to $[L_1(z) \colon \dots \colon L_k(z)]$.
Let us consider the set 
\[
X'=\{
(p,r) \in \P^2 \times \P^{k-1}| L(p)=c_n(r)
\}
\]
and the two restrictions to $X'$ of the projections on the first and second factors of $\P^2 \times \P^{k-1}$, obtaining the following diagram. 
\[
\begin{tikzcd}
X' \arrow[dotted,two heads]{d}{\chi'} \arrow[dotted]{r}{} 
& \P^{k-1} \arrow[,two heads]{d}{c_n}  \\
\P^2 \arrow{r}{L} & \P^{k-1}
\end{tikzcd}
\]
The set $X'$ is defined so that the diagram commutes. 
The map $\chi'$ is a branched covering map of the same degree as $c_n$, ramified along the lines $L_1, \dots, L_k$ in $\P^2$.

Now, one can prove the following. 
\begin{lemma}
A point $q \in X'$ is singular if and only if its image $\chi'(q)$ belongs to at least three lines in the arrangement.
Moreover, the singular points can be resolved by suitable blow-ups.
\end{lemma}

This means that we need to blow up the points of the arrangements where three or more points meet and we will denote $\tau \colon \widehat{\P^2} \to \P^2$ this blow up. 
Moreover, we need to blow up their preimages in $X'$, obtaining a morphism $\rho \colon Y \to X'$. 
Then one can transport the branched covering backwards to $Y$ and obtain the following commuting diagram 
\[
\begin{tikzcd}
Y \arrow[two heads]{d}{\sigma} \arrow{r}{\rho} 
& X' \arrow[,two heads]{d}{\chi'}  \\
\widehat{\P^2} \arrow{r}{L} & \P^2
\end{tikzcd}
\]

Now $\sigma$ is again a branched covering map of the same degree, ramified over the proper transforms in $\widehat{\P^2}$ of the initial lines of the arrangement in $\P^2$ and over the exceptional divisors $\P(\tgt_{p_j}\P^2)$, for each $p_j$ where three lines or more intersect. 
The ramification index at the (proper transform $D_i$ of) each line $L_i$ and at each exceptional divisors $E_j$ (which is the blow up of the point $p_j$) is $n$.

Now take one of the points $p_j$ belonging to $r$ lines of the arrangement ($r \geq 3$) and look at each of the singular points in its preimage. 
When they are resolved with the suitable blow ups, each of them resolves into a smooth curve $C$ and the restriction of $\sigma$ to $C$ is still a branched covering map. 
The degree of the restriction will now be $n^{r-1}$

To summarise, we have the following diagram
\begin{equation}\label{diag}
\begin{tikzcd}
C \arrow[hook]{r}{} \arrow[two heads]{d}{\sigma|_C}&
Y \arrow[two heads]{d}{\sigma} \arrow{r}{\rho} 
& X' \arrow[,two heads]{d}{\chi'} \arrow[]{r}{}
& \P^{k-1} \arrow[two heads]{d}{c_n}
\\
\P(\tgt_p\P^2)\arrow[hook]{r}{} & \widehat{\P^2} \arrow{r} & \P^2 \arrow{r}{L}
& \P^{k-1}
\end{tikzcd}
\end{equation}

In order to verify if a given branched covering is hyperbolic one uses Miyaoka-Yau inequality.  
We state it in the simplest case when the surface is compact. 

\begin{theo}
If $Y$ is a compact surface of general type then $3c_2(Y)-c_1^2(Y)\geq 0$ with equality if and only if $Y$ is complex hyperbolic.
\end{theo}

One can compute the Chern numbers using the combinatorial data of the branched cover.  
When the surface is not compact one has to modify the theorem.

Hirzebruch showed that for $n=5$ and the lines $L$ forming a configuration called the \emph{complete quadrilateral arrangement}, the equality is satisfied.
This means that $Y$ is a complex hyperbolic surface. 
In \cite{YY}, Yamazaki and Yoshida determine the lattice $\Gamma$ in $PU(n,1)$ such that $\widehat{\P^2}$ is the quotient of $\hc$ by $\Gamma$ and the surface $Y$ is the quotient of $\hc$ by the (torsion-free) subgroup $[\Gamma,\Gamma]$. 
In fact, this is a well known lattice, as we will see in Section \ref{sec:DM}.

\begin{figure}[h]
\centering
\includegraphics[width=0.6\textwidth]{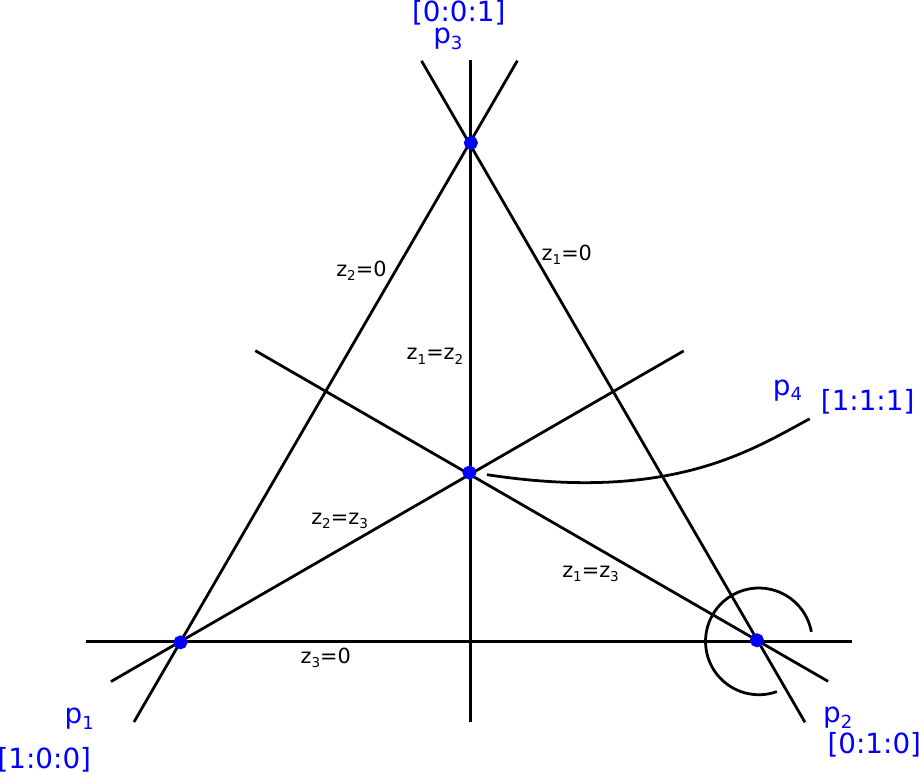}
\begin{quote}\caption{The complete quadrilateral arrangement.\label{fig:complete}} \end{quote}
\end{figure}

The complete quadrilateral arrangement is obtained by taking four points in $\P^2$ in general position and by considering the lines connecting each pair.
This will give six lines with three double intersections and four triple intersection exactly at the initial points.
In some suitable coordinates, we can take points 
\begin{align*}
p_1&=[1\colon 0\colon 0], & p_2&=[0\colon 1 \colon 0], 
& p_3&=[0 \colon 0 \colon 1], & p_4&=[1 \colon 1\colon 1]
\end{align*}
and hence the six lines 
\begin{align*}
l_{ij} &\colon z_i-z_j=0 & &\textit{\emph{ with }} & i,j &\in \{0,1,2,3\}, 
&&i<j
& &\textit{\emph{ and }}& z_0&=0.
\end{align*}
The line arrangement (see Figure \ref{fig:complete}) is then defined by 
\[
z_1z_2z_3(z_2-z_1)(z_3-z_2)(z_1-z_3)=0.
\]

In \cite{BHH} (see also \cite{tretkoff}) Hirzebruch's construction is generalized. 
One defines ramification orders at the complex lines and supposing that a branched cover $Y$ of $\widehat{\P^2}$ exits one looks at ramification orders such that $3c_2(Y)-c_1^2(Y)$ satisfies the correct equality so that $Y$ is a smooth complex hyperbolic surface. One has to prove then the existence of the appropriate branched covering (see \cite{BHH}, \cite{Kobayashi}).

Let $L_1,\cdots, L_k$ be a set of $k$ distinct lines $\C P^2$.   
The intersection points with more than 2 lines are denoted by $p_j$ (called singular points) and the intersection points with only two lines by $p_i$.  
Blowing up the points $p_j$ we obtain a complex surface $X$ with a configuration of divisors (with only normal crossings) which consists of the proper transforms of the lines, noted $D_i$, and the exceptional divisors over the points $p_j$, noted $E_j$.  
For the exceptional divisors we have $E_j\cdot E_j=-1$ and for the proper transforms we have $D_i\cdot D_i=1-\sigma_i$, where $\sigma_i$ is the number of singular points of the arrangement.  
Note then that if $\sigma_i=3$ the proper transform $D_i$ is an exceptional divisor.

One constructs branched coverings $\pi : Y\rightarrow X$ with ramification $m_j$ along $E_j$ and $n_i$ along $D_i$.
 We will always consider branched covers of degree $N$ satisfying the following conditions:
\begin{enumerate}
\item For each point $x\in D_i$ which is not contained in another divisor, $\pi^{-1}(x)$ has $N/n_i$ points. Analogously, for each point $x\in E_j$ which is not contained in another divisor, $\pi^{-1}(x)$ has $N/m_j$ points.
\item For each point $x\in D_i\cap D_j$, $\pi^{-1}(x)$ has $N/n_in_j$ points. Analogously, for each point $x \in D_i\cap E_j$, $\pi^{-1}(x)$ has $N/n_im_j$ points.
\end{enumerate}

The most important arrangement is the complete quadrilateral (also called $\overline{Ceva(1)}$ or $Ceva(2)$ arrangement). Together with the four blown-up points which are intersections of three lines we obtain the arrangements of 10 exceptional divisors.

The ramification index at the exceptional divisors $E_j$ (which is the blow up of the point $p_j$) is set to be the value $m_j$ satisfying the following equation 
\begin{equation}\label{indexBU}
\frac{2}{m_j}+\sum_{i=1}^{r}\frac{1}{n_i}=r-2,
\end{equation}
where $r$ is the number of lines of the arrangement passing through the point $p_j$ and the sum is over these lines (in this case r=3).  This simplifies computations of Chern numbers and one obtains a list of 
branched covers of $\widehat{\P}^2$.
The orders of the ramifications at the ten exceptional divisors are given in Table \ref{table:list} (see \cite{BHH} and \cite{tretkoff}). 

\begin{rk}\label{rk:neginf}
 A negative order at a blown-up point corresponds to a branched cover where, over an exceptional divisor $E$, exceptional divisors appear again.  
To obtain a complex hyperbolic surface one needs to blow down these rational curves and singularities appear.  
The orbifold subgroup preserving $E$ is a finite triangle group. 
On the other hand, if the ramification order at $E$ is infinite, $E$ is covered by elliptic curves.  
Now, the complex hyperbolic surface of finite volume will be the complement of the elliptic curves, the branched covering being a compactification of the complex hyperbolic surface.
 
\end{rk}

\begin{center}
\begin{longtable}{|c||c|c|c|c|c|c|c|c|c|c|}
\hline
$(p,k)$ or $(p,k,p')$ & $n_{*0}$ & $n_{*1}$ & $n_{*2}$ & $n_{*3}$ & $n_{01}$ & $n_{02}$ & $n_{03}$ & $n_{12}$ & $n_{23}$ & $n_{13}$  \\
\hline
$(10,5)$ & 5 & 5 & 5 & 5 & 5 & 5 & 5 & 5 & 5 & 5 \\
$(8,4)$ & 8 & 8 & 8 & 8 & 4 & 4 & 4 & 4 & 4 & 4 \\
(4,8) & -4 & 8 & 8 & 8 & 8 & 8 & 8 & 2 & 2 & 2 \\
(8,2) & 8 & -8 & -8 & -8 & 2 & 2 & 2 & 4 & 4 & 4 \\
(18,3) & 3 & 9 & 9 & 9 & 3 & 3 & 3 & 9 & 9 & 9 \\
(10,2) & 5 & -10 & -10 & -10 & 2 & 2 & 2 & 5 & 5 & 5 \\
(12,4) & 4 & 6 & 6 & 6 & 4 & 4 & 4 & 6 & 6 & 6 \\
(4,4,6) & 12 & 12 & 6 & 6 & 6 & 4 & 4 & 4 & 3 & 4 \\
(12,3) & 4 & 12 & 12 & 12 & 3 & 3 & 3 & 6 & 6 & 6 \\
(12,12,6) & 4 & 4 & 2 & 2 & -6 & 12 & 12 & 12 & 3 & 12 \\
(3,3,4) & -12 & -12 & 12 & 12 & 6 & 3 & 3 & 3 & 2 & 3 \\
(12,2) & 4 & -12 & -12 & -12 & 2 & 2 & 2 & 6 & 6 & 6 \\
(4,6) & -4 & 12 & 12 & 12 & 6 & 6 & 6 & 2 & 2 & 2 \\
(4,3) & -4 & -12 & -12 & -12 & 3 & 3 & 3 & 2 & 2 & 2 \\
(10,3) & 5 & 15 & 15 & 15 & 3 & 3 & 3 & 5 & 5 & 5 \\
(18,2) & 3 & -18 & -18 & -18 & 2 & 2 & 2 & 9 & 9 & 9 \\
(4,5) & -4 & 20 & 20 & 20 & 5 & 5 & 5 & 2 & 2 & 2 \\
(8,3) & 8 & 24 & 24 & 24 & 3 & 3 & 3 & 4 & 4 & 4 \\
(6,3) & $\infty$ & $\infty$ & $\infty$ & $\infty$ & 3 & 3 & 3 & 3 & 3 & 3 \\
(4,4) & -4 & $\infty$ & $\infty$ & $\infty$ & 4 & 4 & 4 & 2 & 2 & 2 \\
(6,6) & $\infty$ & 6 & 6 & 6 & 6 & 6 & 6 & 3 & 3 & 3 \\
(2,6,6) & -6 & 6 & $\infty$ & $\infty$ & 3 & 6 & 6 & 2 & 3 & 2 \\
(6,2) & $\infty$ & -6 & -6 & -6 & 2 & 2 & 2 & 3 & 3 & 3 \\
(6,4) & $\infty$ & 12 & 12 & 12 & 4 & 4 & 4 & 3 & 3 & 3 \\
(3,4,4) & -12 & $\infty$ & 6 & 6 & 12 & 4 & 4 & 3 & 2 & 3 \\
\hline
\caption{Deligne-Mostow lattices and ramification orders.} 
\label{table:list}
\end{longtable}
\end{center}

\begin{figure}[h]
\centering
\includegraphics[width=1\textwidth]{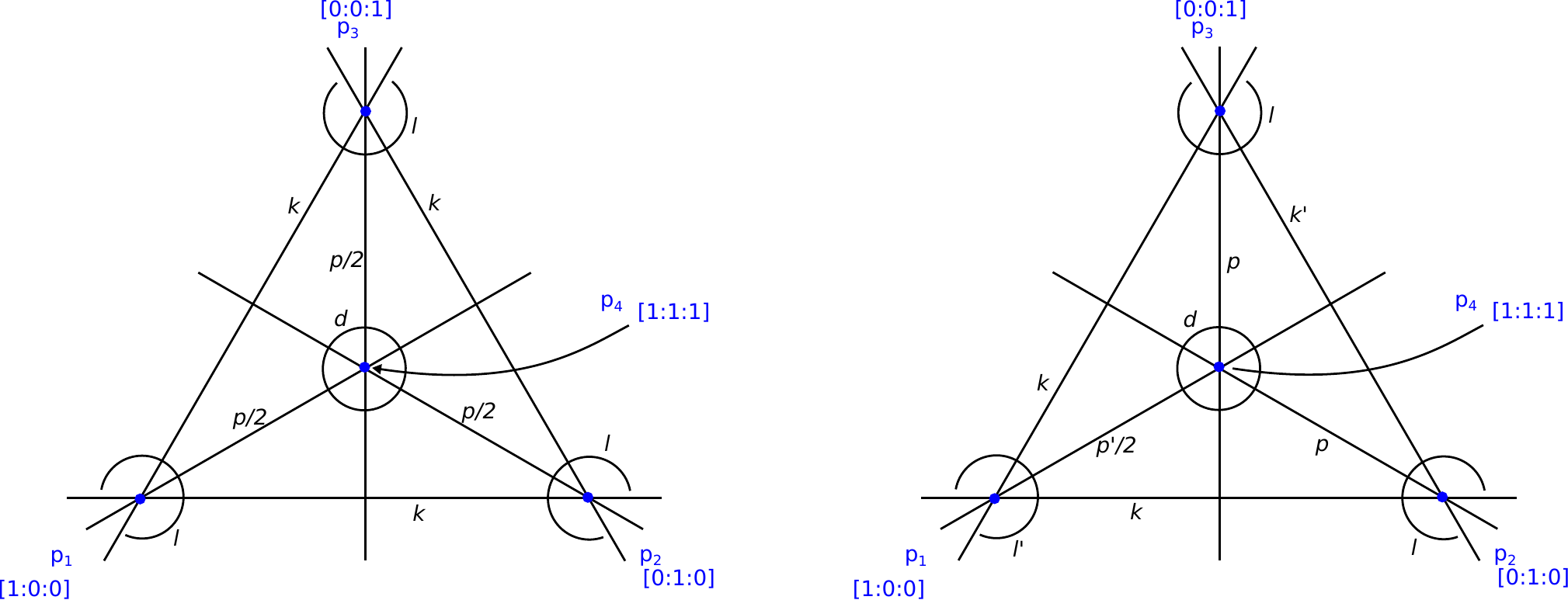}
\begin{quote}\caption{The ramification orders on the line arrangement for the 3- and 2-fold symmetry lattices.\label{fig:arrang}} \end{quote}
\end{figure}

We remark that this list only contains the lattices in Deligne-Mostow list with $P$ even. 
As mentioned by Tretkoff, one could quotient the line arrangement by the symmetry given by the lines with same ramification orders, before taking the branched cover.
The group of symmetries is (a subgroup in the 2-fold symmetry case, of) $ \mathcal S_3$, the permutations on three elements. 
This would allow to take odd values for the ramification orders and hence obtain the full list of Deligne-Mostow lattices.

\subsection{The Deligne-Mostow lattices}\label{sec:DM}

The branched covers $Y$ over the orbifolds $\widehat{\P^2}$ constructed by Tretkoff actually correspond to quotients by (torsion free subgroups of) some well known lattices called \emph{the Deligne-Mostow lattices} in $PU(2,1)$.

Those were introduced by Mostow and then studied by Deligne and Mostow in multiple papers in the 80's using monodromy of hypergeometric functions. 
The 2-dimensional lattices are parametrised by 5 real numbers $\mu=(\mu_0, \dots \mu_4)$ with $0<\mu_i<1$ and $\sum \mu_i=2$. 
They thus obtain a finite list of values of $\mu$ for which one can obtain a lattice.
These lattices are divided in two classes, the lattices with 2- and 3-fold symmetry. 
This means that respectively 2 and 3 of the 5 parameters have the same value. 
The condition on the sum and the symmetry means that the lattices with 2- (resp. 3-) fold symmetry can actually be identified using only 3 (resp. 2) parameters. 

For the 3-fold symmetry case, assume that the three with the same value are $\mu_1=\mu_2=\mu_3$. 
Then we will denote the lattice as $\Gamma_{(p,k)}$, identifying it by the parameters $(p,k)$ chosen in such a way that 
\begin{align*}
\mu_0&=\frac{1}{2}+\frac{1}{p}-\frac{1}{k}, 
& \mu_1=\mu_2=\mu_3&=\frac{1}{2} -\frac{1}{p}, 
& \mu_4=\frac{2}{p}+\frac{1}{k}.
\end{align*}
We will then also be interested in the parameters 
\begin{align}\label{ldef}
\frac{1}{l}&=\frac{1}{2}-\frac{1}{p}-\frac{1}{k}, 
& \frac{1}{d}&=\frac{1}{2} -\frac{3}{p}.
\end{align}

For the 2-fold symmetry case, we assume that the ones with the same values are $\mu_1=\mu_2$. 
Then we will use the triple $(p,k,p')$ to determine the lattice denoted $ \Gamma_{(p,k,p')}$, with the three parameters chosen so that 
\begin{align*}
\mu_0&=\frac{1}{2}+\frac{1}{p'}-\frac{1}{k}, 
& \mu_1=\mu_2&=\frac{1}{2}-\frac{1}{p'}, &
\mu_3&=\frac{1}{2}+\frac{1}{p'}-\frac{1}{p}, 
& \mu_4&=\frac{1}{p}+\frac{1}{k}. \\
\end{align*}
To the lattice $p,k,p'$ we will also associate the parameters
\begin{align*}
\frac{1}{l}&=\frac{1}{2}+\frac{1}{p'}-\frac{1}{p}-\frac{1}{k}, 
& \frac{1}{k'}&=\frac{1}{p}+\frac{1}{k}-\frac{2}{p'}, \\
\frac{1}{l'}&=\frac{1}{2}-\frac{1}{p'}-\frac{1}{k}, 
& \frac{1}{d}&=\frac{1}{2}-\frac{1}{p'}-\frac{1}{p}. \\
\end{align*}

Tretkoff shows the list of the Deligne-Mostow lattices in $PU(2,1)$ is exactly the list that she obtains with the branched covers (see Section \ref{sec:tretkoff}). 
This correspondence is given by using the parameters $p,k,l,d$ for the 3-fold symmetry case and $p,p',k,k',l,l',d$ for the 2-fold symmetry case defined above as the ramification indices as in Figure \ref{fig:arrang}.
In fact, the relation between $\mu$ and the ramification indices is given by 
\begin{equation}\label{TableDM}
\mu_\alpha+\mu_\beta =1- \frac{1}{n_{\alpha\beta}},
\end{equation}
where $n_{\alpha\beta}$ is the weight associated to the line $l_{\alpha\beta}$, for $\alpha, \beta=0,1,2,3$, $\alpha<\beta$.
At the blow up points, we have the associated weight $n_{*\beta}$, with Equation \eqref{TableDM} still holding for $\beta=0,1,2,3$ and defining $\mu_*=\mu_4$ and obtained by $\sum \mu_i=2$.

\subsection{Line arrangements and fundamental domains}\label{fund}

In \cite{thurston}, Thurston reinterpreted these same lattices in terms of moduli spaces of cone metrics on a sphere with fixed cone singularities.
In \cite{irene3F} and \cite{irene2}, an explicit fundamental domain for the lattices is given following his approach. 
Recall that the combinatorial structure of the polyhedra depends on the sign of some parameters. 
Whenever a parameter is negative (resp. infinite), one or more 2-facets collapse to a point (resp. to a point on the boundary).
We now want to show that the line arrangement and the fundamental domains are strictly related. 
We will explain the details in the 3-fold symmetry case, as in the 2-fold symmetry case the construction works in the same way. 

\begin{figure}[h]
\centering
\includegraphics[width=0.6\textwidth]{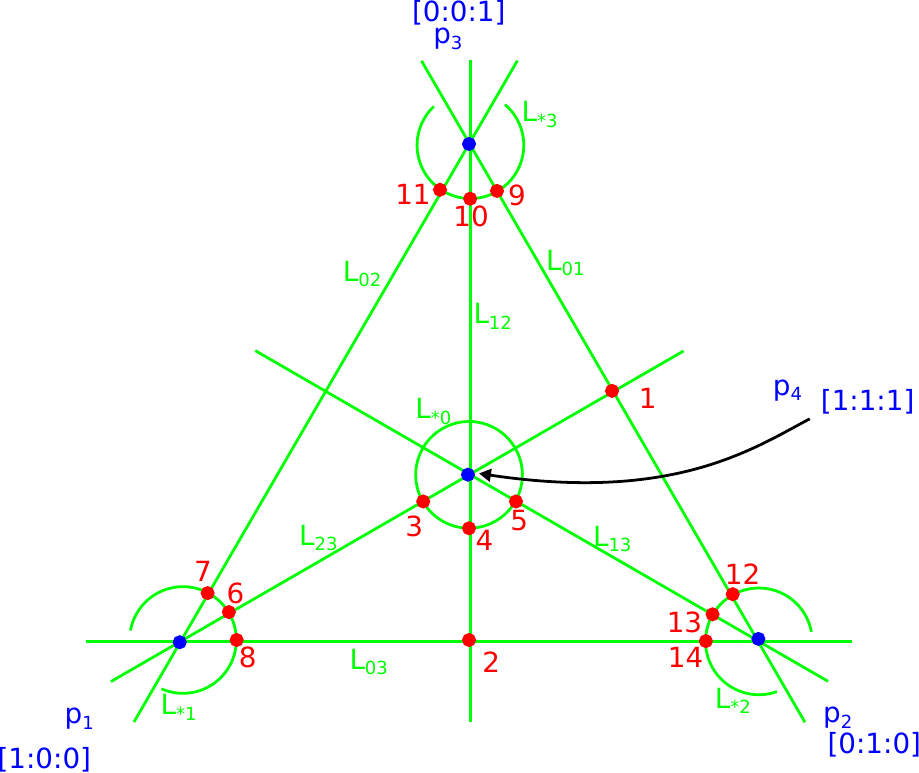}
\begin{quote}\caption{The line arrangement with the corresponding vertices of the polyhedron. \label{LineArrangement}} \end{quote}
\end{figure}

There is a correspondence between the lines in the arrangement $l_{ij}$ defined by $z_i=z_j$ as described above and the configuration of the complex lines $L_{ij}$ where the two cone points $v_i$ and $v_j$ collapse (up to permuting the indices). 
Each of these complex lines contains at most one triangular ridge (a 2-facet of the polyhedron).
The configuration of complex lines is as in Figure \ref{LineArrangement}. 
The numbered dots at the intersection of two lines represents a vertex of the polyhedron. 
Comparing Figure \ref{LineArrangement} and \ref{fig:arrang}

If $L_{ij}$ is the mirror of a certain complex reflection $T$, then the ramification order at the corresponding $l_{ij}$ is the order of $T$.
The triangular ridge on the complex line $L_{ij}$ collapses to a point (potentially on the boundary) exactly when the the ramification order of the corresponding line is negative or infinite. 
These are the cases described in Remark \ref{rk:neginf}.

\begin{rk}\label{rk:orth}
Note that whenever exactly two complex lines in the (blown up) arrangement meet, the two are orthogonal in the fundamental domains picture and if there is a vertex in the intersection, it does not come from collapsing facets.
On the contrary, when a vertex comes from a collapsing ridge, it is at the intersection of three complex lines which meet at an angle, which can be calculated.  
\end{rk}

\begin{figure}[h]
\centering
\includegraphics[width=1\textwidth]{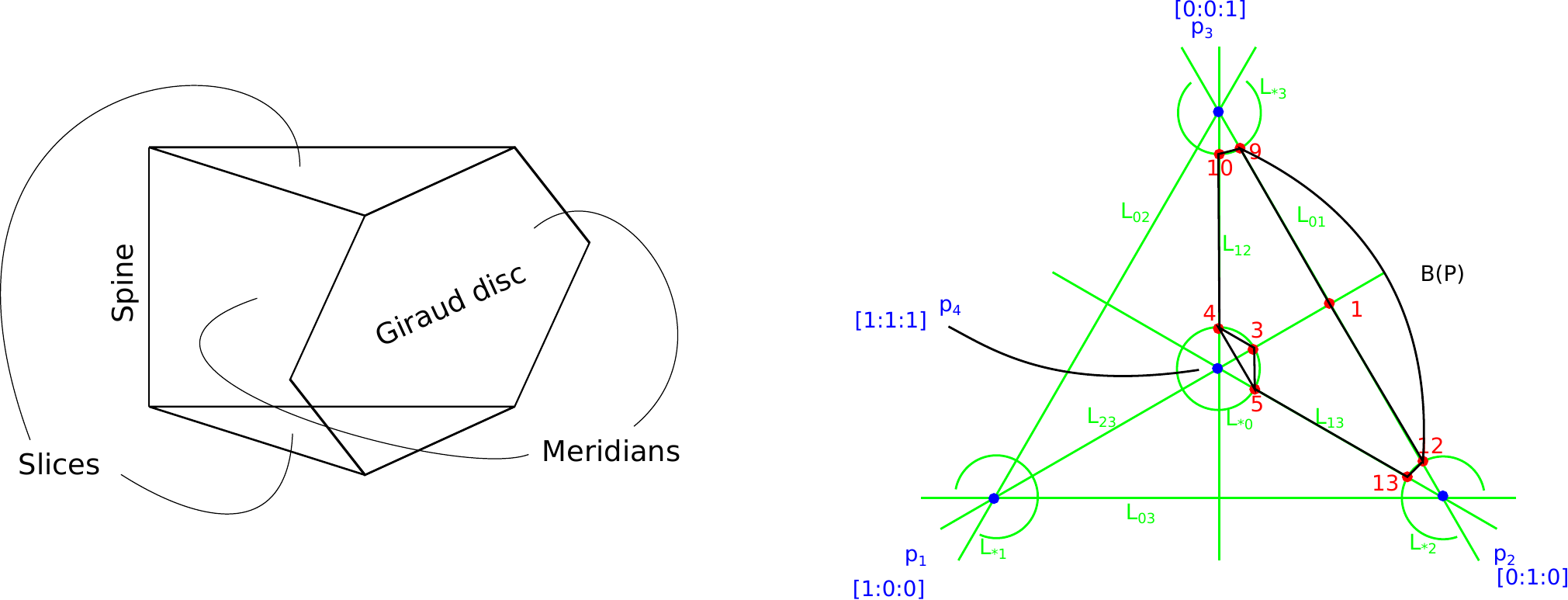}
\begin{quote}\caption{The combinatorial structure of a side and the side $B(P)$, as seen inside the line arrangement. \label{SidesLinesArr}} \end{quote}
\end{figure}

The sides of the polyhedra described in \cite{irene3F} are all contained in bisectors and (before collapsing) all have the same combinatorial structure, which is as in Figure \ref{SidesLinesArr}, where all ridges contained in a side have a special structure. 

By the previous section, we have that $\widehat{\P^2} = \hc/\Gamma$, which is our orbifold (up to forgetting the 3-fold symmetry and taking multiple copies of it). 
When we take a complex line in $\hc$ and we quotient it by $\Gamma$, we obtain a (double) triangle with side identifications that make it equivalent to a copy of $\P^1$, hence a line in $\widehat{\P^2}$. 
In fact, the stabiliser of the complex line is a triangle group (more details will be given in Section \ref{sec:hybrid}).
The information of the ramification orders around the lines record the angles of the initial triangle.

We remark that here we are ignoring the 2- and 3-fold symmetry and taking $R_i^2$ as the complex reflections. 
Since $p$ is the order of $R_i$ and hence $p/2$ is the order of $R_i^2$,this is why Tretkoff only obtains the Deligne-Mostow lattice with $p$ even. 
The complex reflections $T$ at each of the complex line considered are 

\begin{center}
\begin{tabular}{|ccc|ccc|ccc|}
\hline
$L_{ij}$ & $T$ & order & $L_{ij}$ & $T$ & order & $L_{ij}$ & $T$ & order \\
\hline
12 & $R_2$ & $p$ & 01 & $A_1$ & $k$ & *1 & $R_2R_1J$ & $l$ \\
13 & $R_1R_2R_1^{-1}$ & $p$ & 02 & $R_2A_1R_2^{-1}$ & $k$ & *2 & $R_1JR_2$ & $l$ \\
23 & $R_1$ & $p$ & 03 & $R_1R_2A_1R_2^{-1}R_1^{-1}$ & $k$ & *3 & $JR_2R_1$ & $l$ \\
*0 & $P^3$ & $d$ &&&&&& \\
\hline
\end{tabular}
\end{center}

This means that one can also see the sides of the polyhedron inside the line arrangement picture. 
Indeed, the sides of the polyhedron are obtained by considering one of exceptional divisors and a line disjoint from it as top and bottom slices.
An example, using side $B(P)$ in the notation of \cite{irene3F} is given in Figure \ref{SidesLinesArr}.

\section{First application: Representations of 3-manifolds}\label{sec:reps}

\subsection{Lefschetz fibration }

In this section we will introduce a fibration $\widehat{\P^2}\to \P^1$, from which we will deduce a similar fibration $Y \to C$, with $Y$ and $C$ as defined in Section \ref{sec:covers}.
We will then describe the generic and singular fibres, calculating their genus. 
Finally, we will interpret the fibration in terms of the fundamental domains described in Section \ref{fund}. 

\subsubsection{The fibration}

The fibration is described as a \emph{pencil of conics} in $\P^2$ passing through the four points of triple intersection of the arrangement described in Section \ref{sec:covers}. 
It is an example of \emph{Lefschetz fibration}.

A conic is the set of points $[z_1\colon z_2 \colon z_3] \in \P^2$ that are zeros of a quadratic form $q([z_1\colon z_2 \colon z_3])=az_1^2+bz_1z_2+cz_1z_3+dz_2^2+ez_2z_3+fz_3^2$.
The set of all conics is a projective space of dimension 5. 
Imposing the passage through the four points $p_1$, $p_2$, $p_3$ and $p_4$ from Section \ref{sec:covers}%{\color{red}Check that we have the coord of the pi in that section}
 and solving the linear system for the coefficient, one gets a conic of the form 
\begin{equation}\label{eq:conic}
q([z_1\colon z_2 \colon z_3])=c(z_1z_3-z_1z_2)+e(z_2z_3-z_1z_2)=0. 
\end{equation}
Each pair $[c\colon -e]\in \P^1$ corresponds to a single conic passing through the four points. 

%If we take a new point in $\P^2$ different from the previous ones, we can calculate the conic in the set described in \eqref{eq:conic} passing through the fifth point also. 
%We have hence a specific point in $\P^1$. 

We would like to extend this map to the whole of $\P^2$, but at the four intersection points we need more information. 
%We can then specify the tangent direction to the conic at the point $p_i$. 
%This means that instead of associating a conic (hence a point in $\P^1$) to each of the $p_i$,
 We associate a conic to each direction in the tangent space at the point $p_i$, so to each point in $\P(\tgt_{p_i}\P^2)$. 

By putting the two information together, one gets a fibration $f \colon \widehat{\P^2} \to \P^1$ from the blow up at the four points to the set of conics passing through the four points. 
\begin{rk}\label{rk:sections}
The fibration has natural sections $\P^1=\P(\tgt_{p_i}\P^2) \to \widehat{\P^2}$ obtained by identifying $\P^1$ to $\P(\tgt_{p_i}\P^2)$ by associating to each conic the tangent direction at $p_i$ and then see $\P(\tgt_{p_i}\P^2)$ as an exceptional divisor in $\widehat{\P^2}$.
\end{rk}

The fibres of this map are the proper transforms in $\widehat{\P^2}$ of the conics.
Of these, there are exactly three singular ones, each of which is a union of two distinct lines. 
Together, they form the complete quadrilateral arrangement we started from. 

To write the fibration $f$ in coordinates, first let us look at the points in $\P^2$.
One needs to look at Equation \eqref{eq:conic} and remark that in order for the form to be zero, one needs to have $c=-(z_2z_3-z_1z_2)$ and $e=(z_1z_3-z_1z_2)$ and hence we have
\begin{align}\label{eq:fibration}
[z_1\colon z_2 \colon z_3] &\mapsto [z_2(z_1-z_3) \colon z_1(z_2-z_3)].
\end{align}
Remark that this is not defined at the points $p_i$. 
Blowing up at each of these points, one gets coordinate charts $(w_1, w_{2|1})$ and $(w_{1|2}, w_2)$ as given in the Appendix of \cite{ruben}.
For example, blowing up at $p_4$, the fibration on the exceptional divisor is given as as 
\begin{align*}
(w_1, w_{2|1}) &\mapsto [w_1w_{2|1}+1 \colon w_{2|1}(w_1+1) ], \\
(w_{1|2}, w_2) &\mapsto [w_{1|2}(w_2+1) \colon w_{1|2}w_2+1 ].
\end{align*}

Then we have three singular fibres over the points $[0 \colon 1]$, $[0 \colon 0]$ and $[0 \colon 0]$ which are given respectively by $z_2(z_1-z_3)=0$ (hence the union of the two lines $z_2=0$ and $z_1=z_3$), by $z_1(z_2-z_3)=0$ and by $z_3(z_1-z_2)=0$. 

\subsubsection{The fibre}

We now want to apply the fibration $f$ to the blow up of $\P^2$ and compose it with the covering map $\sigma$ from Section \ref{sec:covers}. 
This gives a new fibration $f \circ \sigma \colon Y \to \P^1$. 
Moreover, since $f|_{\P(\tgt_{p_i}\P^2)}$ is an isomorphism, the map $f \circ \sigma|_C$ is a branched covering map with the same properties as $\sigma|_C$.

We hence have the diagram
\begin{equation}\label{fibdiag}
\begin{tikzcd}
C \arrow[hook]{rr}{} \arrow[two heads, swap]{d}{\sigma|_C} \arrow[bend left]{ddrr}{f \circ \sigma|_C}&&
Y \arrow[two heads]{d}{\sigma} \arrow[two heads, bend left]{dd}{f \circ \sigma} 
\\
\P(\tgt_p\P^2)\arrow[hook, crossing over]{rr}{} \arrow{drr}{} && \widehat{\P^2} \arrow[two heads]{d}{f} 
\\
&& \P^1 \arrow{ull}{f|_{\P(\tgt_{p_i}\P^2)}}
\end{tikzcd}
\end{equation}

Recall from Remark \ref{rk:sections} that the inclusion of $\P(\tgt_p\P^2)$ in $\widehat{\P^2}$ is a section of the fibration $f$. 
Similarly, we want to see the inclusion of $C$ in $Y$ as a section of a fibration. 
Indeed, one can prove that there exist a fibration $\tau \colon Y \to C$ which makes the subdiagram of the diagram above only containing $C$, $Y$ and $\P^1$ commute. 
The inclusion of $C$ in $Y$ is a section of $\tau$.

%{\color{red}
%Describe the generic fibre. 
%Calculate its genus. 
%}

\subsubsection{The fundamental domains}

Our goal is to understand the Lefschetz fibration described in \cite{ruben} inside the polyhedron described in \cite{irene3F}. 
The first step to do this is to write the fibration in terms of the $\textbf z$- and $\textbf w$-coordinates that describe the polyhedron. 
We will denote the coordinates in \cite{ruben} in which the fibration is given as $z_i'$ and we will leave the $\textbf z$- and $\textbf w$-coordinates in \cite{irene3F} as they are.

Following Proposition 3.5 of \cite{ruben}, for $n >3$ an integer, we denote $Q_n$ the quotient of $(\P^1)^n$ by the diagonal action of $\aut (\P^1)$ in the sense of geometric invariant theory. 
In other words, one can then think of $Q_5$ as the choice of five points on a sphere, modulo some identifications. 
Then the fibration $f \colon \widehat{\P^2} \to \P^1$ pulls back to a fibration from $Q_5$ to $Q_4$, called the forgetful map, as it ``forgets'' one of the points. 

\begin{equation}\label{diag:forgetful}
\begin{tikzcd}
(v'_1,v'_2,v'_3,v'_4,v'_5) \arrow{rrrr}{} \arrow{ddd}{}&&&& 
\left[ \frac{\det(v'_1,v'_4)}{\det(v'_1,v'_5)} \colon \frac{\det(v'_2,v'_4)}{\det(v'_2,v'_5)} \colon \frac{\det(v'_3,v'_4)}{\det(v'_3,v'_5)} \right] \arrow{ddd}{}
\\
&Q_5 \arrow{rr}{} \arrow{d}{}&& \widehat{\P^2} \arrow{d}{f}&
\\
&Q_4 \arrow{rr}{} && \P^1 &
\\
(v'_1,v'_2,v'_3,v'_4) \arrow{rrrr}{}&&&&
\left[ \frac{\det(v'_1,v'_3)}{\det(v'_1,v'_4)} \colon \frac{\det(v'_2,v'_3)}{\det(v'_2,v'_4)}  \right]
\end{tikzcd}
\end{equation}

To describe the fibration on the $\textbf z$-coordinates, one can use the positions of the points as described in \cite{irene3F}. 
Then one has the expression of the $v_i$ in terms of the $\textbf z$-coordinates.

\begin{align*}
v_*&=0 \\
v_0&=-i  \frac{\sin \theta}{\sin(\theta+\phi)}z_ 3+i  \frac{\sin \alpha}{\sin(\alpha-\phi)}z_ 1, \\
v_1&= -i \frac{\sin \theta}{\sin(\theta+\phi)}z_ 3+i e^{-i\phi} z_ 1, \\
v_2&= -i e^{-i\phi}z_ 2+ie^{-i(\theta+\phi)} \frac{\sin \phi}{\sin(\theta+\phi)}z_ 3, \\
v_3&= -ie^{-i(\theta+\phi)} \frac{\sin \beta}{\sin(\beta-\theta)}z_ 2 +ie^{-i(\theta+\phi)}\frac{\sin \phi}{\sin(\theta+\phi)}z_ 3, \\
\end{align*}

Note that $(v_1', v_2', v_3', v_4', v_5')= (v_1, v_2, v_3, v_0, v_*)$, where, again, the primed expressions correspond to the notation in \cite{ruben} and the other expressions correspond to the notation in \cite{irene3F}. 
The fibration can now be expressed as 
\begin{align}\label{fibration}
f\colon [z'_1:z'_2:z'_3] &\mapsto \left[\frac{v_1-v_3}{v_1-v_0}: \frac{v_2-v_3}{v_2-v_0}\right]
\end{align}
Then from the isomorphism one gets 
\begin{align*}
z_1'&= \frac{v_1-v_0}{v_1-v_*}= 
\frac{\sin(\theta+\phi)\sin\phi(1-e^{i\theta})e^{-i\phi}z_1}
{(\sin\phi+\sin(\theta-\phi))
(\sin(\theta+\phi)e^{-i\phi}z_1-\sin\theta z_3)}, \\
z_2'&= \frac{v_2-v_0}{v_2-v_*}= 
\frac{\sin(\theta+\phi)
(-\sin\theta z_1
-(\sin\phi+\sin(\theta-\phi))e^{-i\phi}z_2
+(\sin\phi+\sin(\theta-\phi))e^{-i\phi}z_3)}
{(\sin\phi+\sin(\theta-\phi))
(-\sin(\theta+\phi)e^{-i\phi}z_2+\sin\phi e^{-i(\theta+\phi)}z_3)}, \\
z_3'&= \frac{v_3-v_0}{v_3-v_*}= 
\frac{\sin(\theta+\phi)
(\sin\theta z_1
-(\sin\phi+\sin(\theta-\phi))e^{-i(\theta+\phi)}z_2
+(\sin\phi+\sin(\theta-\phi))e^{-i\phi}z_3)}
{(\sin\phi+\sin(\theta-\phi))
(-\sin(\theta+\phi)e^{-i(\theta+\phi)}z_2+\sin\phi e^{-i(\theta+\phi)}z_3)}, \\
\end{align*}

Substituting these expressions in the formula for the fibration and remembering the relation between the $\textbf z$- and $\textbf w$-coordinates as in \cite{irene3F}, one gets:

\begin{align*}
f\left(
\begin{bmatrix}
z_1 \\ z_2 \\z_3
\end{bmatrix}
\right)
&= \frac{
(-\sin\theta e^{i\phi}z_1-(\sin\phi+\sin(\theta-\phi))z_2 +(\sin\phi+\sin(\theta-\phi))z_3)
(z_1+e^{-i\theta}z_2-z_3)
}{
-\sin\phi (1-e^{i\theta})(1-e^{-i\theta})z_1z_2
}, \\
f\left(
\begin{bmatrix}
w_1 \\ w_2 \\w_3
\end{bmatrix}
\right)
&= \frac{
(-\sin\theta e^{-i\phi}w_1-e^{-i\theta}(\sin\phi+\sin(\theta-\phi))w_2 +(\sin\phi+\sin(\theta-\phi))w_3)
(w_1+w_2-w_3)
}{
(-\sin\theta e^{-i\phi}w_1-(\sin\phi+\sin(\theta-\phi))w_2 +(\sin\phi+\sin(\theta-\phi))w_3)
(w_1+e^{-i\theta}w_2-w_3)
}.
\end{align*}

Then we have 
\begin{prop}
Inside the fundamental polyhedra $D$, the fibre $[\alpha,\colon 1]\in\P^1 \setminus\{[1 \colon 0]\}$ by the Lefschetz fibration is the solution of the following equation 
\[
f\left(
\begin{bmatrix}
z_1 \\ z_2 \\z_3
\end{bmatrix}
\right)=\alpha,
\]
for $f$ as above.
\end{prop}

We know that the three singular fibres are the preimages of the points $[1:0]$, $[0:1]$ and $[1:1] \in \P^1$. 
When writing $f$ as in \eqref{fibration}, it is immediate to see that the singular fibres are as follows: 
\begin{align}\label{singularfibres}
f^{-1}([1:0]) &=L_{23}\cup L_{01}, 
& f^{-1}([0:1]) &=L_{13}\cup L_{02}, 
& f^{-1}([1:1]) &=L_{12}\cup L_{03}.
\end{align}

Using this, we can look at the projection of the boundary of the polyhedron. 
The vertices of the polyhedron all belong to one of the complex lines that form the singular fibres as in \eqref{singularfibres}. 
More precisely, 
\begin{align*}
&\textbf z_1, \textbf z_3, \textbf z_6, \textbf z_9, \textbf z_{12}, 
& &\textit{ \emph{projects to} } & [1:0], \\
&\textbf z_2, \textbf z_4, \textbf z_8, \textbf z_{10}, \textbf z_{14}, 
& &\textit{ \emph{projects to} } & [1:1], \\
&\textbf z_5, \textbf z_7, \textbf z_{11}, \textbf z_{13}, 
& &\textit{ \emph{projects to} } & [0:1], \\
\end{align*}

Similarly, of all the edges in the polyhedron, %(see Table page 89 of thesis), 
only few do not project down on single points of $\P^1$.
They project instead to a curve between two of the marked points $0, 1, \infty \in \P^1$ as recorded in the following table.

\begin{center}
\begin{tabular}{|cc|cc|cc|}
\hline
Curve & Edges & Curve & Edges & Curve & Edges \\
\hline
$1 - \infty$ & $\gamma_{3,4}$ 
& $0 - 1$ & $\gamma_{4,5}$
& $0 - \infty$ & $\gamma_{3,5}$ \\
& $\gamma_{6,8}$
&& $\gamma_{7,8}$
&& $\gamma_{6,7}$ \\
& $\gamma_{9,10}$
&& $\gamma_{10,11}$
&& $\gamma_{9,11}$ \\
& $\gamma_{12,14}$
&& $\gamma_{13,14}$
&& $\gamma_{12,13}$ \\
\hline
\end{tabular}
\end{center}

One can then use the formulae for $f$ and calculate the actual curve in $\P^1$ that the edges project to. 
For example, consider a point $\textbf z \in \gamma_{3,4}$. 
Then 
\begin{align*}
\textbf z
&=\begin{bmatrix}
z_1 \\ z_2 \\z_3
\end{bmatrix}=
\begin{bmatrix}
\frac{\sin\phi+\sin(\theta-\phi)}{\sin(\theta+\phi)}\\ u \\1
\end{bmatrix},
& &\textit{\emph{where}}
& u \in \left(
0, \frac{\sin\phi(2\cos\theta-1)}{\sin(\theta+\phi)}
\right)
\end{align*}

and it projects to the curve

\[
\frac{(\sin\phi e^{-i\theta}-\sin(\theta+\phi)u)
(\sin\phi(1-2\cos\theta)+\sin(\theta+\phi)e^{-i\theta}u)}
{-\sin\phi \sin(\theta+\phi)(1-e^{-i\theta})(1-e^{i\theta})u}.
\]

Note that when $u=0$, then $\textbf z$ is the vertex $\textbf{z}_3$ and it projects to $\infty$. 
When $u=\frac{\sin\phi(2\cos\theta-1)}{\sin(\theta+\phi)}$, then $\textbf z$ is the vertex $\textbf{z}_4$ and it projects to $1$. 
This indeed projects to a curve between 1 and $\infty$.

The same can be done for each other edge.
Then we get: 
\begin{prop}
The 2-facets of the polyhedron contained in complex lines through the Lefschetz fibration project either to one of the marked points $[1 \colon 0], [0 \colon 1], [1 \colon 1]\in \P^1$ to triangles in $\P^1$ of vertices the same marked points.
\end{prop}

\subsubsection{Representation of 3-manifolds}\label{sub:representations}

In this section we review Dashyan's construction of representations of 3-manifolds fundamental groups in $PU(2,1)$ and obtain a description of further representations associated to Deligne-Mostow lattices.

Denote a Lefshetz fibration $\tau : Y\rightarrow C$ and define $Y^u\subset Y$ to be the complement of the singular fibres.  
We now have a fibration $\tau_{\vert Y^u} Y^u \rightarrow C^u$ with no singular fibres.  
Consider now an element $[\gamma] \in \pi_1(C^u)$, where $\gamma : S^1 \rightarrow C^u$ is a loop.  
Let $M_\gamma$ be the the 3-manifold obtained by pulling-back the fibre bundle $Y$ through the loop $\gamma$.  
In the particular case that $\gamma$ is an embedding, $M_\gamma$ will be embedded in $Y$. Denote by $F_0$ a fixed non-singular fibre so that $\pi_1(M_\gamma)=\left<\gamma\right>\ltimes\pi_1(F_0)$.

If $Y$ is a complex hyperbolic manifold, we can write $\rho : \pi_1(Y)\rightarrow PU(2,1)$ to denote a particular embedding of the the fundamental group so that $Y=\hc/\rho(\pi_1(Y))$.  
Composing with the natural map $\pi_1(M_\gamma)\rightarrow \pi_1(Y)$, we obtain a representation 
$$
\rho_\gamma : \pi_1(M_\gamma)\rightarrow PU(2,1).
$$

Consider the case of the Lefshetz fibration
$\widehat{\P^2}\rightarrow {\P^1}$ whose generic fibre is a sphere with four marked points. Let $\Mod_{0,4}$ be the mapping class group of the sphere with four marked points.

\begin{prop} [\cite{ruben}] Let
$\widehat{\P^2}^u\rightarrow {\P^1}^u$ be the fibration with no singular fibres.  
Then the monodromy representation $\pi_1({\P^1}^u)\rightarrow \Mod_{0,4}$ is an isomorphism.
\end{prop}

\begin{prop}[\cite{ruben}] 
For any element $[\gamma]$ in $\pi_1(C^u)$, if its image in $\pi_1(C)$ is not trivial, then
\begin{enumerate}
\item the kernel of $\rho_\gamma$ is equal to the kernel of $\pi_1(F_0)\to\pi_1(Y)$,
\item the monodromy of the fibration ${Y_1}^u\to C^u$ along $\gamma$ is pseudo-Anosov,
\item the kernel is not of finite type.
\end{enumerate}
\end{prop}

For $Y$ obtained as a branched covering described by the Ceva(2) configuration, one has the following description. The proof follows the same lines as in \cite{ruben} once we know that the quotient of hyperbolic space by a Deligne-Mostow lattice is an orbifold described by a Lefschetz fibration $\widehat{\P^2}\rightarrow {\P^1}$ with regular fibres which are sphere-orbifolds with four marked points with isotropy with order equal to the multiplicity of the branching order at the exceptional divisors.

\begin{theo}
For each element of $\Mod_{0,4}$, consider a surface bundle $M$ with that monodromy and with fibre the orbifold with the sphere as underlying space and with isotropy of order equal to the multiplicity of the exceptional fibre at each of the four marked points. Then there exists a representation of the orbifold fundamental group of $M$ into a Deligne-Mostow lattice in $\Isom(\hc)$.
\end{theo}

\section{Second application: complex hybridisation}\label{sec:hybrid}

One of the main problems in complex hyperbolic geometry is the existence of non-arithmetic lattices. 
In fact, few examples are known and only in dimension 2 and 3. 
In dimension 4 or higher, existence of non-arithmetic lattices is still an open question.
The hope for the existence comes from the parallel with the real hyperbolic space, in which infinitely many non-arithmetic lattices exist in any dimension. 
This has been proved in \cite{nonarithm} by Gromov and Piatetski-Shapiro using a technique called \emph{hybridisation}. 
It consists in taking two arithmetic lattices of dimension $n$ and produce a new lattice by ``glueing'' along a hypersurface. 
When the two lattices are non-commensurable, the resulting lattice is shown to be non-arithmetic. 

Since in complex hyperbolic space there are no totally geodesic real hypersurfaces, an immediate analogue of the hybridisation construction is not possible. 
Nonetheless, possible analogues have been explored first by Paupert (which he attributes to Hunt) and then by Paupert and Wells. 
Here we will use Wells' definition from \cite{wells}, which is the following. 

\begin{defin}
Let $\Gamma_1, \Gamma_2 <PU(n,1)$ be lattices. 
The \emph{hybrid} of $\Gamma_1$ and $\Gamma_2$ is a group $H=H(\Gamma_1,\Gamma_2)=\langle \Lambda_1, \Lambda_2 \rangle$, where $\Lambda_1, \Lambda_2 <PU(n+1,1)$ are two discrete subgroups which stabilise totally geodesic (complex) hypersurfaces $\Sigma_1$ and $\Sigma_2$ which are orthogonal. Moreover, we require that $\Gamma_i=\Lambda_i |_{\Sigma_i}$ and that the intersection $\Lambda_1 \cap \Lambda_2< PU(n-1,1)$ is a lattice.
\end{defin}

In his paper, Wells proves that two of the non-arithmetic Deligne-Mostow lattices can be given as hybrids (in this sense) of two non-commensurable arithmetic lattices in $PU(1,1)$.

We now want to use the relation between Hizebruch's construction and the fundamental domains in \cite{irene3F} to prove that the procedure from Wells can be extended to more of the 2-dimensional Deligne-Mostow lattices. 
In fact, Wells did not have explicit fundamental domains available and so there was the extra difficulty of finding suitable orthogonal subspaces along which to hybridise. 
The fundamental domains suggest which pairs of complex lines are orthogonal and using Hirzebruch's construction it is immediate to see which are the triangle group associated to them. 
We will use the parameters $(p,k)$ and $(p,k,p')$ as in Section \ref{sec:DM} to determine the lattices. 

There are 9 commensurability classes in the 2-dimensional Deligne-Mostow lattices. 
Of these, two have been treated by Wells. 
They are $(4,6)$ (commensurable to $(4,4,6)$), and $(5,4)$. 
Here we will show that more of them can be given as hybrids. 
Moreover, no other of the Deligne-Mostow lattices can be give as a hybrid of the groups relates to a line of our arrangement.
More precisely, we have that 
\begin{theo}\label{theo:hybrid}
Let $\Gamma\subset PU(2,1)$ be a non-arithmetic Deligne-Mostow lattice.
Let $(L_1, L_2)$ be an orthogonally intersecting pair of lines with corresponding triangle groups $\Delta_1$ and $\Delta_2$.  
Then 
\begin{enumerate}
\item 
 either both $\Delta_1$ and $\Delta_2$ are arithmetic and commensurable,
 \item or both $\Delta_1$ and $\Delta_2$ are arithmetic and non-commensurable and in this case $\Gamma$ is an hybridisation of $\Delta_1$ and $\Delta_2$,
 \item or one of the triangle groups is non-arithmetic.
\end{enumerate}

\end{theo}

The second case corresponds to the hybridisation construction and we can prove exactly in which cases this happens.   

\begin{corollary}
Other than the two cases $(4,6)$ and $(5,4)$ found by Wells, also the lattices $(3,8)$, $(6,4)$ and $(3,4,4)$ are hybrids of two non-commensurable arithmetic lattices in $PU(1,1)$. 
\end{corollary}

%{\color{red}Is the fact that the commensurable ones have multiple ways of hybridising intereresting?}

\subsection{Wells' cases}

In \cite{wells}, Wells looks at subspaces $e_i^\perp$ and $v_{ijk}^\perp$, which are the mirrors of $R_i$ and $J^{\pm}R_jR_k$, with $k \equiv i \pm 1 \textit{ \emph{(mod 3)}}$, respectively. 
Note that his notation for the generators is standard and it is the same used in \cite{irene3F}. 
This means that we can use the presentation given with the fundamental domain to determine which ridges (2-facets, hence which complex lines in the arrangement) correspond to Wells' subspaces.

He chooses subspaces $v_{312}^\perp$ and $v_{321}^\perp$ as orthogonal complex lines along which to hybridise. 
They first one is the mirror of $J^{-1}R_1R_2=J^{-1}P$, which is the cycle transformation of the ridge $F(P,J)$, which is contained in $L_{01}$. 
The second one is the mirror of $JR_2R_1$, which is the cycle transformation associated to $F(R_1,J^{-1})$, which is contained in $L_{*3}$. 

Looking at Figure \ref{fig:arrang}, first we remark that the values considered are only for $p<6$, so $d<0$ and hence $L_{*0}$ is blown down to a single point. 
Moreover, $L_{*3}$ and $L_{01}$ are indeed orthogonal by Remark \ref{rk:orth}.
Moreover, if we look at the stabiliser of a line in the arrangement, we are looking at the triangle group $\Delta(i,j,k)$, where $i, j, k$ are the ramification orders of the three complex lines intersecting the initial line. 
Remark that when we have a triangle group with two equal angles, we can reduce the triangle to its half. 
This means to consider another triangle group commensurable to the first one. 
Then the order at the apex of the triangle doubles and from $\Delta(i,i,j)$ we get $\Delta(2,2i,j)$.
So, for example, the stabiliser of $v_{312}^\perp$ and hence $L_{01}$ will be $\Delta(2,p,l)$ because $L_{01}$ intersects $L_{*3}$ and $L_{*2}$, which both have ramification order $l$, and also $L_{23}$, which has ramification order $p/2$. 
Similarly, the stabiliser of $v_{321}^\perp$ and hence $L_{*3}$ will be $\Delta(2,p,k)$. 
These are indeed the values given by Wells in Tables 2, 3 and 4.

\subsection{Proof of Theorem \ref{theo:hybrid}}

In order to see which case of the theorem holds, we need to choose pairs of orthogonal lines and calculate the triangle groups in $PU(1,1)$ giving their stabilisers in the group. 
We will include only the calculations for the second case (which also prove the corollary). 
For the others, one uses the same procedure to see which triangle groups are obtained and then compare them to the table page 418 of \cite{trgps}, which gives arithmeticity and commensurability classes of triangle groups. 
%Then we will show that the group generated by these two is the full lattice. 

\subsubsection*{3-fold symmetry cases}

Like in Wells' cases, let us look at the orthogonal complex lines $L_{*3}$ and $L_{01}$. 
Note that the lattices $(3,8)$ and $(6,4)$ both belong to the lattices of second type, so the lines $L_{*i}$ for $i=1,2,3$ do not collapse (i.e. blow down), while the line $L_{*0}$ collapses to a point inside $\hc$ in the first case and on the boundary in the second case. 
This means that we can always look at the two lines $L_{*3}$ and $L_{01}$. 

As in the previous section, the subgroups preserving line $L_{ij}$ is generated by the complex reflections whose mirrors intersect $L_{ij}$. 
Moreover, two of these are enough, because the other will be a combination of these. 
So the stabiliser of $L_{*3}$ is generated by, for example, $J^{-1}P$ and $R_2$, which have $L_{01}$ and $L_{12}$ as mirrors, respectively. 
The stabiliser of $L_{01}$ is now generated by $JR_2R_1$ and $R_1$, which have $L_{*3}$ and $L_{23}$ as mirrors, respectively. 
For $(3,8)$, they are triangle groups $\Delta(2,3,24)$ and $\Delta(2,3,8)$, which are both arithmetic and non-commensurable (cfr table page 418 of \cite{trgps}).
For $(6,4)$ they are $\Delta(2,6,12)$ and $\Delta(2,4,6)$, which are again both arithmetic and non-commensurable.

Let us hence look at the hybrid group of these two $\Gamma =\langle J^{-1}P,R_2,JR_2R_1, R_1 \rangle < \Gamma_{(p,k)}$.
From the presentation of $\Gamma_{(p,k)}$ given in \cite{irene3F}, we can see that $R_1$, $R_2$ and $J$ are enough to generate the full group. 
But now $J= J^{-1}PR_2^{-1}R_1^{-1}$ and hence the two groups are the same. 

This means that the two non-arithmetic Deligne-Mostow lattices with 3-fold symmetry in $PU(2,1)$, $(3,8)$ and $(6,4)$ can be given as hybrids of two non-commensurable arithmetic lattices in $PU(1,1)$. 

\subsubsection*{2-fold symmetry case}

In the 2-fold symmetry case we have more choices of orthogonal lines, since the lines arrangement has fewer symmetries. 
Here we will choose $L_{03}$ and $L_{02}$, although multiple other choices would give similar hybridisation results. 

The lines intersecting $L_{01}$ are $L_{*3}$, $L_{23}$ and $L_{*2}$, which have ramification orders $l$, $p'/2$ and $l$ respectively. 
The stabiliser of $L_{01}$ is hence a triangle group $\Delta(2,p',l)$. 
Now $L_{*3}$ is the mirror of $R_0K$ and $L_{23}$ is the mirror of $R_0$. 

Similarly, the lines intersecting $L_{*3}$ are $L_{01}$, $L_{12}$ and $L_{02}$, which have ramification orders $k$, $p$ and $k'$ respectively. 
The stabiliser of $L_{*3}$ is hence a triangle group $\Delta(k,p,k')$. 
Now $L_{01}$ is the mirror of $Q^{-1}K=A_1$ and $L_{12}$ is the mirror of $B_2$. 

Let us now look at $\Gamma=\langle R_0K, R_0, A_1, B_2 \rangle <\Gamma_{(p,k,p')}$.
We want to show that $\Gamma$ is the full $\Gamma_{(p,k,p')}$. 
Now, using the presentation in \cite{irene2}, we can see that $\Gamma_{(p,k,p')}$ is generated by $B_1$, $R_0$ and $A_1$, so all we need is to express $B_1$ in terms of the generators of $\Gamma_{(p,k,p')}$. 
But from the presentation we can see that $B_1=R_0B_2R_0^{-1}$ and hence we are done. 

Now for the lattice $(4,4,6)$, we have that $\Delta(2,p',l)=\Delta(2,6,6)$ and $\Delta(k,p,k')=\Delta(4,4,6)$ are both arithmetic and non-commensurable. 
Similarly, for the lattice $(3,4,4)$, we have that $\Delta(2,p',l)=\Delta(2,4,6)$ and $\Delta(k,p,k')=\Delta(4,3,12)$ are both arithmetic and non-commensurable. 
This proves the theorem. 

Note that $(4,4,6)$ is commensurable for $(4,6)$ so this just gives another way to see it as a hybrid.

\newpage
\addcontentsline{toc}{section}{\refname}
\bibliographystyle{alpha}
\bibliography{biblio}

\begin{thebibliography}{BHH87}

\bibitem[BHH87]{BHH}
Gottfried Barthel, Friedrich Hirzebruch, and Thomas H\"{o}fer.
\newblock {\em Geradenkonfigurationen und {A}lgebraische {F}l\"{a}chen}.
\newblock Aspects of Mathematics, D4. Friedr. Vieweg \& Sohn, Braunschweig,
  1987.

\bibitem[Das19]{ruben}
Ruben Dashyan.
\newblock A construction of representations of 3-manifold groups into pu(2,1)
  through lefschetz fibrations.
\newblock {\em arXiv preprint arXiv:1909.10229}, 2019.

\bibitem[DM86]{delignemostow}
P.~Deligne and G.~D. Mostow.
\newblock Monodromy of hypergeometric functions and nonlattice integral
  monodromy.
\newblock {\em Inst. Hautes \'Etudes Sci. Publ. Math.}, (63):5--89, 1986.

\bibitem[GPS88]{nonarithm}
M.~Gromov and I.~Piatetski-Shapiro.
\newblock Nonarithmetic groups in {L}obachevsky spaces.
\newblock {\em Inst. Hautes \'{E}tudes Sci. Publ. Math.}, (66):93--103, 1988.

\bibitem[Hir83]{hirz}
F.~Hirzebruch.
\newblock Arrangements of lines and algebraic surfaces.
\newblock In {\em Arithmetic and geometry, {V}ol. {II}}, volume~36 of {\em
  Progr. Math.}, pages 113--140. Birkh\"{a}user, Boston, Mass., 1983.

\bibitem[Kob90]{Kobayashi}
Ryoichi Kobayashi.
\newblock Uniformization of complex surfaces.
\newblock In {\em K\"{a}hler metric and moduli spaces}, volume~18 of {\em Adv.
  Stud. Pure Math.}, pages 313--394. Academic Press, Boston, MA, 1990.

\bibitem[Mos80]{mostow}
George~D. Mostow.
\newblock On a remarkable class of polyhedra in complex hyperbolic space.
\newblock {\em Pacific J. Math.}, 86(1):171--276, 1980.

\bibitem[MR03]{trgps}
Colin Maclachlan and Alan~W. Reid.
\newblock {\em The arithmetic of hyperbolic 3-manifolds}, volume 219 of {\em
  Graduate Texts in Mathematics}.
\newblock Springer-Verlag, New York, 2003.

\bibitem[Par06]{livne}
J.~R. Parker.
\newblock Cone metrics on the sphere and {L}ivn\'e's lattices.
\newblock {\em Acta Math.}, 196(1):1--64, 2006.

\bibitem[Par09]{survey}
John~R. Parker.
\newblock Complex hyperbolic lattices.
\newblock In {\em Discrete groups and geometric structures}, volume 501 of {\em
  Contemp. Math.}, pages 1--42. Amer. Math. Soc., Providence, RI, 2009.

\bibitem[Pas16]{irene3F}
Irene Pasquinelli.
\newblock Deligne-{M}ostow lattices with three fold symmetry and cone metrics
  on the sphere.
\newblock {\em Conform. Geom. Dyn.}, 20:235--281, 2016.

\bibitem[Pas19]{irene2}
Irene Pasquinelli.
\newblock Fundamental domains and presentations for the {D}eligne--{M}ostow
  lattices with 2-fold symmetry.
\newblock {\em Pacific J. Math.}, 302(1):201--247, 2019.

\bibitem[Thu98]{thurston}
William~P. Thurston.
\newblock Shapes of polyhedra and triangulations of the sphere.
\newblock In {\em The {E}pstein birthday schrift}, volume~1 of {\em Geom.
  Topol. Monogr.}, pages 511--549. Geom. Topol. Publ., Coventry, 1998.

\bibitem[Tre16]{tretkoff}
Paula Tretkoff.
\newblock {\em Complex ball quotients and line arrangements in the projective
  plane}, volume~51 of {\em Mathematical Notes}.
\newblock Princeton University Press, Princeton, NJ, 2016.
\newblock With an appendix by Hans-Christoph Im Hof.

\bibitem[Wel19]{wells}
Joseph Wells.
\newblock Non-arithmetic hybrid lattices in pu(2,1).
\newblock {\em arXiv preprint arXiv:1905.12184}, 2019.

\bibitem[YY84]{YY}
Tadashi Yamazaki and Masaaki Yoshida.
\newblock On {H}irzebruch's examples of surfaces with {$c^{2}_{1}=3c_{2}$}.
\newblock {\em Math. Ann.}, 266(4):421--431, 1984.

\end{thebibliography}
\nocite{*}

\end{document}